# RANDOM ORIENTED TREES: A MODEL OF DRAINAGE NETWORKS

By Sreela Gangopadhyay, Rahul Roy[1] and Anish Sarkar[2]

*Indian Statistical Institute*

Consider the $d$-dimensional lattice $\mathbb{Z}^d$ where each vertex is "open" or "closed" with probability $p$ or $1-p$, respectively. An open vertex $v$ is connected by an edge to the closest open vertex $w$ such that the $d$th co-ordinates of $v$ and $w$ satisfy $w(d) = v(d) - 1$. In case of nonuniqueness of such a vertex $w$, we choose any one of the closest vertices with equal probability and independently of the other random mechanisms. It is shown that this random graph is a tree almost surely for $d = 2$ and 3 and it is an infinite collection of distinct trees for $d \geq 4$. In addition, for any dimension, we show that there is no bi-infinite path in the tree and we also obtain central limit theorems of (a) the number of vertices of a fixed degree $\nu$ and (b) the number of edges of a fixed length $l$.

**1. Introduction.** Leopold and Langbein (1962) introduced a geometric model of natural drainage network which they described as

> *using a sheet of rectangular cross-section graph paper, each square is presumed to represent a unit area. Each square is to be drained, but the drainage channel from each square has equal chance of leading off in any of the four cardinal directions, subject only to the condition that, having made a choice, flow in the reverse direction is not possible. Under these conditions it is possible for one or more streams to flow into a unit area, but only one can flow out.*

Subsequently Scheidegger (1967) introduced a direction of flow. In his study of Alpine valleys, he imposed conditions on the Leopold and Langbein model by requiring that the drainage paths be in the "direction of high gradients between watershed and main valleys." Thus the drainage forms an oriented

Received April 2002; revised July 2003.
[1]Supported in part by DST Grant MS/152/01.
[2]Supported by Indian National Science Academy, Brazilian Science Academy and IME, Universidade de São Paulo.

network, with a square emptying to one of its two neighbors in a preferred direction. Howard (1971) removed the restriction of drainage to a neighboring square and modelled a network to include "headward growth and branching in a random fashion." Rodriguez-Iturbe and Rinaldo (1997) present a survey of the development of this field.

The random graph we study here follows the one described by Howard (1971) with the caveat that a stream is not permitted to terminate or become inactive. Thus we consider the $d$-dimensional lattice $\mathbb{Z}^d$ where each vertex is "open" or "closed" with probability $p$ or $1-p$, respectively. The open vertices represent the water sources. An open vertex $v$ is connected by an edge to the closest open vertex $w$ such that the $d$th co-ordinates of $v$ and $w$ satisfy $w(d) = v(d) - 1$. In case of nonuniqueness of such a vertex $w$, we choose any one of the closest vertices with equal probability and independently of the other random mechanisms. These edges represent the channels of flow in the drainage network.

Our main result (Theorem 2.1) is that, for $d = 2$ and 3, all the tributaries connect to form one single delta, while for $d \geq 4$, there are infinitely many deltas, each with its own distinct set of tributaries. In this connection it is worth noting that (Theorem 2.2) there is no main river, in the sense that there is no bi-infinite river; instead, each tributary has its own distinct source. In addition, for any dimension, we obtain central limit theorems of (a) the number of sites where a fixed number $\nu$ of tributaries drain, as well as of (b) the number of channels of a fixed length $l$.

Similar tree–forest dichotomies have been studied for the uniform spanning tree model by Pemantle (1991) and for the minimal spanning tree model by Newman and Stein (1996). Ferrari, Landim and Thorisson (2002) have obtained similar results for a continuous version of this model.

In two dimensions we obtain the main result by showing that the distance between two streams starting at two different sites forms a martingale and thereby invoking the martingale convergence theorem. For three dimensions we employ a technique based on Lyapunov functions, while in four or higher dimensions we couple the streams starting at two different sites with two independent and identically distributed random walks starting at these two sites. To show that there are no bi-infinite paths in the graph we utilize the stationarity of the model and use a Burton–Keane type argument. The limit theorems are obtained by checking that the random processes satisfy the conditions needed to apply Lyapunov's central limit theorem.

The formal details of the model and the statements of results are in the next section.

## 2. The model and statement of results.

Let $\Omega = \{0,1\}^{\mathbb{Z}^d}$ and let $\mathcal{F}$ be the $\sigma$ algebra generated by finite-dimensional cylinder sets. On $(\Omega, \mathcal{F})$ we assign a product probability measure $P_p$ which is defined by its marginals as

$$P_p\{\omega:\omega(u)=1\}=1-P_p\{\omega:\omega(u)=0\}=p \qquad \text{for } u\in\mathbb{Z}^d \text{ and } 0\le p\le 1.$$

Let $\{U_{u,v}:u,v\in\mathbb{Z}^d, v(d)=u(d)-1\}$ be i.i.d. uniform $(0,1]$ random variables on some probability space $(\Xi,\mathcal{S},\mu)$. Here and subsequently we express the co-ordinates of a vector $u$ as $u=(u(1),\dots,u(d))$.

Consider the product space $(\Omega\times\Xi,\mathcal{F}\times\mathcal{S},\mathbb{P}:=P_p\times\mu)$. For $(\omega,\xi)\in\Omega\times\Xi$, let $\mathcal{V}(=\mathcal{V}(\omega,\xi))$ be the random vertex set defined by

$$\mathcal{V}(\omega,\xi)=\{u\in\mathbb{Z}^d:\omega(u)=1\}.$$

Note that if $u\in\mathcal{V}(\omega,\xi)$ for some $\xi\in\Xi$, then $u\in\mathcal{V}(\omega,\xi')$ for all $\xi'\in\Xi$ and thus we say that a vertex $u$ is open in a configuration $\omega$ if $u\in\mathcal{V}(\omega,\xi)$ for some $\xi\in\Xi$.

For $u\in\mathbb{Z}^d$, let

$$\begin{aligned} N_u &= N_u(\omega,\xi)\\ &=\Bigg\{v\in\mathcal{V}(\omega,\xi):v(d)=u(d)-1 \text{ and}\\ &\qquad \sum_{i=1}^d|v(i)-u(i)|=\min\Bigg\{\sum_{i=1}^d|w(i)-u(i)|:w\in\mathcal{V}(\omega,\xi),\\ &\qquad\qquad w(d)=u(d)-1\Bigg\}\Bigg\}.\end{aligned}$$

Note that for $p>0$, $N_u$ is nonempty almost surely and that $N_u$ is defined for all $u$, irrespective of it being open or closed. For $u\in\mathbb{Z}^d$, let

$$\begin{aligned} &h(u)\in N_u(\omega,\xi) \quad \text{be such that}\\ &U_{u,h(u)}(\xi)=\min\{U_{u,v}(\xi):v\in N_u(\omega,\xi)\}.\end{aligned} \tag{1}$$

Again note that for $p>0$ and for each $u\in\mathbb{Z}^d$, $h(u)$ is open, almost surely unique and $h(u)(d)=u(d)-1$. On $\mathcal{V}(\omega,\xi)$ we assign the edge set $\mathcal{E}=\mathcal{E}(\omega,\xi):=\{\langle u,h(u)\rangle:u\in\mathcal{V}(\omega,\xi)\}$.

Consider that graph $\mathcal{G}=(\mathcal{V},\mathcal{E})$ consisting of the vertex set $\mathcal{V}$ and edge set $\mathcal{E}$. For $p=0$, $\mathcal{V}=\varnothing$ almost surely, and, for $p=1$, $\langle u,v\rangle\in\mathcal{E}$ if and only if $u(i)=v(i)$ for all $i\ne d$ and $|u(d)-v(d)|=1$. Also, for a vertex $u\in\mathcal{V}(\omega,\xi)$, there is exactly one edge "going down" from $u$; that is, there is a unique edge $\langle u,v\rangle$ with $v(d)\le u(d)$. Thus the graph $\mathcal{G}$ contains no loops almost surely. Hence, for $0<p<1$, the graph $\mathcal{G}$ consists of only trees. Our first result is THEOREM 2.1. *Let $0 < p < 1$. For $d = 2$ and $d = 3$, $\mathcal{G}$ consists of one single tree $\mathbb{P}$-almost surely, while for $d \geq 4$, $\mathcal{G}$ is a forest consisting of infinitely many disjoint trees $\mathbb{P}$-almost surely.*

Regarding the geometric structure of the graph $\mathcal{G}$, we have

THEOREM 2.2. *Let $0 < p < 1$. For any $d \geq 2$, the graph $\mathcal{G}$ contains no bi-infinite path $\mathbb{P}$-almost surely.*

Now for $\nu \geq 0$, let $S_n$ be the number of vertices in $\mathcal{V} \cap ([1, n]^d)$ of the graph $\mathcal{G}$ with degree $\nu + 1$. Also, for $l \geq 1$, let $L_n$ be the number of edges of $L_1$-length $l$ in the graph $\mathcal{G}$ with one end vertex in $\mathcal{V} \cap ([1, n]^d)$.

THEOREM 2.3. *As $n \to \infty$:*

(a) $\frac{S_n - E(S_n)}{n^{d/2}}$ *converges weakly to a normal random variable;*
(b) $\frac{L_n - E(L_n)}{n^{d/2}}$ *converges weakly to a normal random variable.*

Finally, for $d = 2$, given that a vertex $v$ is open, the following proposition gives the exact distribution of the degree of $v$.

PROPOSITION 2.1. *Given that a vertex $v$ is open, the degree of the vertex in the graph $\mathcal{G}$ has the same distribution as that of $1 + Y + X_1 + X_2$, where $Y$, $X_1$ and $X_2$ are independent nonnegative random variables such that*

$$Y = \begin{cases} 0, & \text{with probability } 1 - p, \\ 1, & \text{with probability } p, \end{cases}$$

$$\mathbb{P}(X_1 \geq r) = \mathbb{P}(X_2 \geq r) = \begin{cases} 1, & \text{for } r = 0, \\ \dfrac{(1-p)^{2r-1}(2-p)}{2(3 - 3p + p^2)^r}, & \text{for } r \geq 1. \end{cases}$$

*Thus the expected degree of a vertex, given that it is open, is 2.*

REMARK 2.1. As in Lemma 7 of Aldous and Steele (1992), using the ergodicity of the process, it may be shown that in any dimension, the expected degree of a vertex, given that it is open, is 2.

**3. Proof of Theorem 2.1.** We fix $0 < p < 1$ and for $u, v \in \mathbb{Z}^{d-1}$ consider the $d$-dimensional vectors $\mathbf{u} := (u, 0)$ and $\mathbf{v} := (v, 0)$ and let $(X_u^n, -n) := h^n(\mathbf{u})$, where $h^n$ denotes the $n$-fold composition of $h$ defined in (1). For $Z_n (= Z_n(u, v)) := X_u^n - X_v^n$, we first observe that it is a time-homogeneous Markov chain with state space $\mathbb{Z}^{d-1}$; indeed, this follows on writing $\{Z_{n+1} = z_{n+1}, Z_n = z_n, \ldots, Z_0 = z_0\} = \bigcup_{x_{n+1}, \ldots, x_0 \in \mathbb{Z}^d} \{X_{x_0}^{n+1} = x_{n+1}, X_{x_0+z_0}^{n+1} = x_{n+1} + z_{n+1},$ $X^n_{x_0} = x_n, X^n_{x_0+z_0} = x_n + z_n, \dots, X^0_{x_0} = x_0, X^0_{x_0+z_0} = x_0 + z_0\}$ and using the Markovian property of the process $\{(X^n_u, X^n_v) : n \geq 0\}$.

The connectedness or otherwise of the graph $\mathcal{G}$ is equivalent to whether or not $Z_n$ is absorbed at the origin. For $d = 2$ and 3, we show that $Z_n$ gets absorbed at the origin, $\mathbf{0} \in \mathbb{Z}^{d-1}$ with probability 1; while for $d \geq 4$, $Z_n$ is a transient Markov chain and hence has a positive probability of not being absorbed. In this connection observe that instead of the above $Z_n$, if we had considered a modified Markov chain $\tilde{Z}_n$, where $\mathbf{0}$ is no longer an absorbing state, but from $\mathbf{0}$ we move in one step to some fixed vertex $\mathbf{u} \neq \mathbf{0}$ with probability 1 and the other transition probabilities are kept unchanged, then to show that the original process $Z_n$ is absorbed at $\mathbf{0}$ almost surely, it suffices to show that the modified Markov process $\tilde{Z}_n$ is recurrent. A more formal argument for this would require $Z_n$ and $\tilde{Z}_n$ to be coupled together until they hit the origin, which occurs almost surely if the modified process is recurrent. For the case $d = 3$, we will show that $\tilde{Z}_n$ is recurrent. The proof is divided into three sections according as $d = 2$, $d = 3$ and $d \geq 4$.

3.1. $d = 2$. Fix $i < j$ and observe that $X^n_i \leq X^n_j$ for every $n \geq 1$, where $X^n_i$ and $X^n_j$ are as defined earlier. Thus the Markov chain $Z_n := X^n_j - X^n_i$ with $Z_0 = j - i$ has as its state space the set of all nonnegative integers. Since the marginal distributions of the increments of $X^n_i$ and $X^n_j$ are identical with finite means, $\{Z_n : n \geq 0\}$ is a nonnegative martingale. Hence, by the martingale convergence theorem [see Billingsley (1979), Theorem 35.4, page 416], $Z_n$ converges almost surely as $n \to \infty$. Since $\{Z_n : n \geq 0\}$ is also a time-homogeneous Markov chain with 0 as the only absorbing state, we must have $Z_n \to 0$ as $n \to \infty$ with probability 1. Since this is true for all $i < j$, we have the result for $d = 2$.

3.2. $d = 3$. Throughout this section the letters $\mathbf{u}$, $\mathbf{v}$ in bold font denote vectors in $\mathbb{Z}^3$, u, v in roman font denote vectors in $\mathbb{Z}^2$ and $u$, $v$ in italic font denote integers. Fix two vectors $\mathbf{u} := (\mathrm{u}, 0)$ and $\mathbf{v} := (\mathrm{v}, 0)$ in $\mathbb{Z}^2 \times \{0\}$ and let $\tilde{Z}_n (= \tilde{Z}_n(\mathrm{u}, \mathrm{v}))$ be the time-homogeneous Markov chain with state space $\mathbb{Z}^2$ as defined at the beginning of this section. We shall exhibit, by a Lyapunov function technique, that this Markov chain $\tilde{Z}_n$ is recurrent, thereby showing that $Z_n$ is absorbed at the origin with probability 1.

Consider the function $f : \mathbb{R}^2 \to [0, \infty)$ defined by $f(\mathrm{x}) := \sqrt{\log(1 + \|\mathrm{x}\|_2^2)}$ where $\| \cdot \|_2$ is the standard $L_2$ norm (Euclidean distance). Since $f(\mathrm{x}) \to \infty$ as $\|\mathrm{x}\|_2 \to \infty$, by Foster's criterion [see Asmussen (1987), Proposition 5.3 of Chapter I, page 18] the following lemma implies that $\tilde{Z}_n$ is recurrent.

LEMMA 3.1. *For all $n \geq 0$, there exists $T \geq 0$ such that, for all $\|\mathrm{x}\|_2 \geq T$, we have*

$$\mathbb{E}(f(\tilde{Z}_{n+1}) - f(\tilde{Z}_n) \mid \tilde{Z}_n = \mathrm{x}) < 0.$$

PROOF. Let $g\colon[0,\infty)\to[0,\infty)$ be defined as $g(x):=\sqrt{\log(1+x)}$. Clearly $g(x)\geq 0$ for all $x\geq 0$ and $g(x)\to\infty$ as $x\to\infty$. Also, for $x,y\geq 0$, the Taylor series expansion yields

$$g(x)-g(y)\leq (x-y)g^{(1)}(y)+\frac{(x-y)^2}{2}g^{(2)}(y)+\frac{(x-y)^3}{6}g^{(3)}(y), \tag{2}$$

which holds because the fourth derivative

$$\begin{aligned} g^{(4)}(s)&=-\frac{3}{(1+s)^4g(s)}-\frac{11}{4(1+s)^4(g(s))^3}\\ &\quad-\frac{18}{8(1+s)^4(g(s))^5}-\frac{15}{16(1+s)^4(g(s))^7}<0 \qquad \text{for } s>0. \end{aligned}$$

The first three derivatives of $g$, which we will be needing shortly, are

$$\begin{aligned} g^{(1)}(s)&=\frac{1}{2(1+s)g(s)},\\ g^{(2)}(s)&=-\frac{1}{2(1+s)^2g(s)}-\frac{1}{4(1+s)^2(g(s))^3},\\ g^{(3)}(s)&=\frac{1}{(1+s)^3g(s)}+\frac{3}{4(1+s)^3(g(s))^3}+\frac{3}{8(1+s)^3(g(s))^5}. \end{aligned}$$

Note that, for all $s$ large,

$$g^{(3)}(s)\leq\frac{3}{(1+s)^3g(s)}.$$

Assuming for the moment that (we will prove this shortly), for some $\alpha>0$,

$$\mathbb{E}(\|\tilde{Z}_{n+1}\|_2^2-\|\tilde{Z}_n\|_2^2\mid\tilde{Z}_n=\mathrm{x})=\alpha+o(\|\mathrm{x}\|_2^{-2}), \tag{3}$$

$$\mathbb{E}((\|\tilde{Z}_{n+1}\|_2^2-\|\tilde{Z}_n\|_2^2)^2\mid\tilde{Z}_n=\mathrm{x})\geq 2\alpha\|\mathrm{x}\|_2^2, \tag{4}$$

$$\mathbb{E}((\|\tilde{Z}_{n+1}\|_2^2-\|\tilde{Z}_n\|_2^2)^3\mid\tilde{Z}_n=\mathrm{x})=O(\|\mathrm{x}\|_2^2), \tag{5}$$

as $\|\mathrm{x}\|_2\to\infty$, and using the above estimates and expression for derivatives, we have, for all $\beta:=\|\mathrm{x}\|_2^2$ large and for some nonnegative constants $C_1$ and $C_2$,

$$\begin{aligned} &\mathbb{E}(f(\tilde{Z}_{n+1})-f(\tilde{Z}_n)|\tilde{Z}_n=\mathrm{x})\\ &\quad\leq\frac{\alpha+C_1/\beta}{2(1+\beta)\sqrt{\log(1+\beta)}}-\frac{2\alpha\beta}{4(1+\beta)^2\sqrt{\log(1+\beta)}}\\ &\qquad-\frac{2\alpha\beta}{8(1+\beta)^2\sqrt{(\log(1+\beta))^3}}+\frac{3C_2\beta}{(1+\beta)^3\sqrt{\log(1+\beta)}}\\ &\quad=\frac{1}{8(1+\beta)^2\sqrt{\log(1+\beta)}}\left[4\alpha+4C_1+\frac{4C_1}{\beta}+\frac{24C_2\beta}{1+\beta}-\frac{2\alpha\beta}{\log(1+\beta)}\right]. \end{aligned}$$

The term inside the square brackets tends to $-\infty$ as $\beta \to \infty$; therefore, for all sufficiently large $\beta$, the term is negative. Thus to complete the proof of the lemma we need to show (3)–(5).

Let $D_k := \{\mathrm{v} \in \mathbb{Z}^2 : \|\mathrm{v}\|_1 \leq k\}$ denote "$L_1$-diamond" of radius $k$ and $\delta D_k := \{\mathrm{v} \in \mathbb{Z}^2 : \|\mathrm{v}\|_1 = k\}$ its boundary, where $\|\cdot\|_1$ denotes the $L_1$ norm. Consider the probability distribution of the step size of the random walk, associated with the tree generated by one particle, that is, the distribution of $X_{\mathrm{o}}^1$:

$$\begin{aligned} p_{\mathrm{u}} &:= \mathbb{P}(X_{\mathrm{o}}^1 = \mathrm{u}) \\ &= \begin{cases} p, & \text{if } \mathrm{u} = \mathrm{o}, \\ \dfrac{(1-p)^{\#D_{k-1}}(1-(1-p)^{\#\delta D_k})}{\#\delta D_k}, & \text{for } \mathrm{u} \in \delta D_k,\ k \geq 1, \end{cases} \end{aligned} \tag{6}$$

where $\mathrm{o} := (0,0)$ is the origin and $\#A$ denotes the cardinality of the set $A$. For any $k \geq 1$ and $i, j \geq 0$, define

$$m_i(k) := \sum_{\mathrm{u}:=(u_1,u_2)\in D_k} u_1^i \, p_{\mathrm{u}}$$

and

$$m_{i,j}(k) := \sum_{\mathrm{u}:=(u_1,u_2)\in D_k} u_1^i u_2^j \, p_{\mathrm{u}}.$$

Since $(-u_1, -u_2) \in D_k$ whenever $(u_1, u_2) \in D_k$, it is clear that, for every $k \geq 1$ we have

$$\begin{aligned} m_i(k) &= 0 \qquad \text{for all odd } i \quad \text{and} \\ m_{i,j}(k) &= 0 \qquad \text{whenever either } i \text{ or } j \text{ is odd.} \end{aligned} \tag{7}$$

Further, since $\#D_k = 1 + 2k(k+1)$ and $\#\delta D_k = 4k$, we have that, for all even $i$,

$$\begin{aligned} 0 < m_i := \lim_k m_i(k) &= \sum_{\mathrm{u}:(u_1,u_2)\in\mathbb{Z}^2} u_1^i p_{\mathrm{u}} \\ &\leq \sum_{k=1}^{\infty} (\#\delta D_k)(\max\{u_i : (u_1,u_2) \in D_k\})^i \, \frac{(1-p)^{\#D_{k-1}}(1-(1-p)^{\#\delta D_k})}{\#\delta D_k} \\ &= \sum_{k=1}^{\infty} k^i (1-p)^{1+2k(k-1)}(1-(1-p)^{4k}) \\ &< \infty. \end{aligned}$$

Similarly, when both $i$ and $j$ are even, $m_{i,j}(k) \to m_{i,j}$ as $k \to \infty$, where $0 < m_{i,j} := \sum_{\mathrm{u}\in\mathbb{Z}^2} u_1^i u_2^j p_{\mathrm{u}} < \infty$. Further, $p_{\mathrm{u}}$ being the same for every $\mathrm{u}$ on $\delta D_k$, the various quantities $m_i$ and $m_{i,j}$ remain unchanged if in their definitions we had considered $u_2$ instead of $u_1$.

Moreover,

$$k^2(m_2 - m_2(k)) \le k^2 \sum_{u:(u_1,u_2)\notin D_k} u_1^2 p_u$$
$$\le \sum_{j=k+1}^{\infty} j^4 (1-p)^{1+2j(j-1)} \to 0$$

as $k \to \infty$ since the sum

$$\sum_{j=1}^{\infty} j^4(1-p)^{1+2j(j-1)} < \infty.$$

A similar result holds for $m_0(k)$ and so we have

$$m_2(k) = m_2 + o(k^{-2}) \quad \text{and} \quad m_0(k) = m_0 + o(k^{-2}) \qquad \text{as } k \to \infty. \tag{8}$$

Now we proceed to compute the expectations:

$$\begin{aligned} &\mathbb{E}(\|\tilde{Z}_{n+1}\|_2^2 - \|\tilde{Z}_n\|_2^2 \mid \tilde{Z}_n = \mathrm{x}) \\ &\quad = \sum_{\mathrm{a},\mathrm{b}\in\mathbb{Z}^2} (\|\mathrm{x}+\mathrm{a}-\mathrm{b}\|_2^2 - \|\mathrm{x}\|_2^2) \\ &\qquad\quad \times \sum_{\mathrm{w}\in\mathbb{Z}^2} [\mathbb{P}\{X_{\mathrm{u}}^{n+1} = X_{\mathrm{u}}^n + \mathrm{a}, \\ &\qquad\qquad\qquad X_{\mathrm{v}}^{n+1} = X_{\mathrm{v}}^n + \mathrm{b} \mid X_{\mathrm{v}}^n = \mathrm{w}, X_{\mathrm{u}}^n = \mathrm{w}+\mathrm{x}\}] \\ &\qquad\qquad \times [\mathbb{P}\{X_{\mathrm{v}}^n = \mathrm{w}, X_{\mathrm{u}}^n = \mathrm{w}+\mathrm{x} \mid \tilde{Z}_n = \mathrm{x}\}] \\ &\quad = \sum_{\mathrm{a},\mathrm{b}\in\mathbb{Z}^2} (\|\mathrm{x}+\mathrm{a}-\mathrm{b}\|_2^2 - \|\mathrm{x}\|_2^2) \\ &\qquad\quad \times \mathbb{P}\{X_{\mathrm{x}}^1 = \mathrm{x}+\mathrm{a}, X_{\mathrm{o}}^1 = \mathrm{b} \mid X_{\mathrm{o}}^0 = \mathrm{o}, X_{\mathrm{x}}^0 = \mathrm{x}\}, \end{aligned} \tag{9}$$

where we have used the translation invariance of the model.

To calculate the above sum we let $k := \|\mathrm{x}\|_2/4$. Note, for $\mathrm{a},\mathrm{b} \in D_k$, we have $\mathbb{P}\{X_{\mathrm{x}}^1 = \mathrm{x}+\mathrm{a}, X_{\mathrm{o}}^1 = \mathrm{b} \mid X_{\mathrm{o}}^0 = \mathrm{o}, X_{\mathrm{x}}^0 = \mathrm{x}\}0 = p_{\mathrm{a}}p_{\mathrm{b}}$; thus, using (7) and (8),

$$\begin{aligned} T_1(1) := &\sum_{\mathrm{a},\mathrm{b}\in D_k} (\|\mathrm{x}+\mathrm{a}-\mathrm{b}\|_2^2 - \|\mathrm{x}\|_2^2) \\ &\qquad \times \mathbb{P}\{X_{\mathrm{x}}^1 = \mathrm{x}+\mathrm{a}, X_{\mathrm{o}}^1 = \mathrm{b} \mid X_{\mathrm{o}}^0 = \mathrm{o}, X_{\mathrm{x}}^0 = \mathrm{x}\} \end{aligned}$$

(10)

$$
\begin{aligned}
&= \sum_{\mathrm{a},\mathrm{b}\in D_k} [(a_1-b_1)^2 + 2x_1(a_1-b_1) \\
&\qquad\qquad + (a_2-b_2)^2 + 2x_2(a_2-b_2)]p_{\mathrm{a}}p_{\mathrm{b}} \\
&= 4m_2(k)m_0(k) \\
&= 4m_2 + o(k^{-2}) \qquad \text{as } k\to\infty.
\end{aligned}
$$

Also, if $\mathrm{b}\notin D_k$, then, taking $\|\mathrm{b}\|_1 = k+l$ for some $l\geq 1$, the occurrence of the event $\{X_{\mathrm{o}}^1=\mathrm{b}\}$ requires that all the vertices in the diamond $D_{k+l-1}$ be closed and that at least one vertex of $\delta D_{k+l}$ be open—an event which occurs, with probability $(1-p)^{1+2(k+l-1)(k+l)} - (1-p)^{1+2(k+l)(k+l+1)}$. Moreover, if $\{X_{\mathrm{o}}^1=\mathrm{b}\}$ occurs, then $X_{\mathrm{x}}^1$ must lie in the smallest diamond centered at x which contains the vertex b; thus $\|X_{\mathrm{x}}^1 - X_{\mathrm{o}}^1\|_2 \leq \|X_{\mathrm{x}}^1\|_1 + \|X_{\mathrm{o}}^1\|_1 \leq (\|\mathrm{x}\|_1 + \|\mathrm{b}\|_1) + \|\mathrm{b}\|_1 = 6k+2l$. Now noting that there are $4(k+l)$ vertices on $\delta D_{k+l}$ and that an argument similar to the above may be given when $\mathrm{a}\notin D_k$, we have

$$
\begin{aligned}
T_2(1) &:= \sum_{\mathrm{a}\notin D_k \text{ or } \mathrm{b}\notin D_k} (\|\mathrm{x}+\mathrm{a}-\mathrm{b}\|_2^2 - \|\mathrm{x}\|_2^2) \\
&\qquad\qquad \times \mathbb{P}\{X_{\mathrm{x}}^1 = \mathrm{x}+\mathrm{a}, X_{\mathrm{o}}^1 = \mathrm{b} \mid X_{\mathrm{o}}^0 = \mathrm{o}, X_{\mathrm{x}}^0 = \mathrm{x}\} \\
&\leq 2\sum_{l\geq 1} 4(k+l)((6k+2l)^2 + (4k)^2)(1-p)^{1+2(k+l-1)(k+l)} \\
&\qquad \times [1-(1-p)^{4(k+l)}] = o(k^{-2}) \qquad \text{as } k\to\infty.
\end{aligned} \tag{11}
$$

This establishes (3) with $\alpha = 4m_2$.

For (4), calculations as in (9) show that $\mathbb{E}((\|\tilde{Z}_{n+1}\|_2^2 - \|\tilde{Z}_n\|_2^2)^2 \mid \tilde{Z}_n = \mathrm{x}) \geq T_1(2)$ where, performing calculations as in (10),

$$
\begin{aligned}
T_1(2) &:= \sum_{\mathrm{a},\mathrm{b}\in D_k} (\|\mathrm{x}+\mathrm{a}-\mathrm{b}\|_2^2 - \|\mathrm{x}\|_2^2)^2 \mathbb{P}\{X_{\mathrm{x}}^1 = \mathrm{x}+\mathrm{a}, X_{\mathrm{o}}^1 = \mathrm{b} \mid X_{\mathrm{o}}^0 = \mathrm{o}, X_{\mathrm{x}}^0 = \mathrm{x}\} \\
&= \sum_{\mathrm{a},\mathrm{b}\in D_k} [(a_1-b_1)^2 + 2x_1(a_1-b_1) + (a_2-b_2)^2 + 2x_2(a_2-b_2)]^2 p_{\mathrm{a}}p_{\mathrm{b}} \\
&= 8m_2(k)m_0(k)\|x\|_2^2 + 4m_4(k)m_0(k) + 16(m_2(k))^2 + 4m_{2,2}(k)m_0(k).
\end{aligned}
$$

Now, as $k\to\infty$, from (8) we have $m_2(k)m_0(k)\|x\|_2^2 - m_2\|x\|_2^2 = 16(m_2 + o(k^{-2}))(1+o(k^{-2}))k^2 - 16m_2k^2 = 16k^2o(k^{-2}) + o(k^{-2}) = o(1)$; also $4m_4(k)m_0(k) + 16(m_2(k))^2 + 4m_{2,2}(k)m_0(k) \to 4m_4 + 16m_2^2 + 4m_{2,2} > 0$. This establishes (4).

Finally, for (5), we write $\mathbb{E}((\|\tilde{Z}_{n+1}\|_2^2 - \|\tilde{Z}_n\|_2^2)^3 \mid \tilde{Z}_n = \mathrm{x}) = T_1(3) + T_2(3)$, where

$$\begin{aligned}
T_1(3) &:= \sum_{\mathrm{a},\mathrm{b}\in D_k} (\|\mathrm{x}+\mathrm{a}-\mathrm{b}\|_2^2 - \|\mathrm{x}\|_2^2)^3 \mathbb{P}\{X_{\mathrm{x}}^1 = \mathrm{x}+\mathrm{a}, X_{\mathrm{o}}^1 = \mathrm{b} \mid X_{\mathrm{o}}^0 = 0, X_{\mathrm{x}}^0 = \mathrm{x}\} \\
&= \sum_{\mathrm{a},\mathrm{b}\in D_k} [(a_1-b_1)^2 + 2x_1(a_1-b_1) + (a_2-b_2)^2 + 2x_2(a_2-b_2)]^3 p_{\mathrm{a}} p_{\mathrm{b}} \\
&= \sum_{\mathrm{a},\mathrm{b}\in D_k} p_{\mathrm{a}} p_{\mathrm{b}} [(a_1-b_1)^3(2x_1+a_1-b_1)^3 \\
&\qquad + 3(a_1-b_1)^2(2x_1+a_1-b_1)^2(2x_2+a_2-b_2)(a_2-b_2) \\
&\qquad + 3(a_1-b_1)(2x_1+a_1-b_1)(2x_2+a_2-b_2)^2(a_2-b_2)^2 \\
&\qquad + (a_2-b_2)^3(2x_2+a_2-b_2)^3] \\
&= T_{1,1}(3) + T_{1,2}(3) + T_{1,3}(3) + T_{1,4}(3) \qquad \text{(say)}.
\end{aligned}$$

Now,

$$\begin{aligned}
T_{1,1}(3) &:= \sum_{\mathrm{a},\mathrm{b}\in D_k} p_{\mathrm{a}} p_{\mathrm{b}} (a_1-b_1)^3 [(a_1-b_1)^3 + 6x_1(a_1-b_1)^2 \\
&\qquad + 12x_1^2(a_1-b_1) + 8x_1^3] \\
&= \sum_{\mathrm{a},\mathrm{b}\in D_k} p_{\mathrm{a}} p_{\mathrm{b}} (a_1-b_1)^6 + 12x_1^2 \sum_{\mathrm{a},\mathrm{b}\in D_k} p_{\mathrm{a}} p_{\mathrm{b}} (a_1-b_1)^4 \\
&= C_3(k) + C_4(k) x_1^2,
\end{aligned}$$

where $C_3(k)$ and $C_4(k)$ are both polynomials in $m_i(k)$ and $m_{i,j}(k)$, each of which converges to the corresponding polynomial in $m_i$ and $m_{i,j}$ as $k \to \infty$. Similar calculations show that $T_{1,4}(3) := \sum_{\mathrm{a},\mathrm{b}\in D_k} p_{\mathrm{a}} p_{\mathrm{b}} (a_2-b_2)^3[(a_2-b_2)^3 + 6x_2(a_2-b_2)^2 + 12x_2^2(a_2-b_2) + 8x_2^3] = C_3(k) + C_4(k)x_2^2$.

Also,

$$\begin{aligned}
T_{1,2}(3) &= 3 \sum_{\mathrm{a},\mathrm{b}\in D_k} p_{\mathrm{a}} p_{\mathrm{b}} (a_1-b_1)^2 [(a_1-b_1)^2 + 4x_1(a_1-b_1) + 4x_1^2] \\
&\qquad \times [2x_2(a_2-b_2) + (a_2-b_2)^2] \\
&= 3 \sum_{\mathrm{a},\mathrm{b}\in D_k} p_{\mathrm{a}} p_{\mathrm{b}} (a_1-b_1)^4 (a_2-b_2)^2 \\
&\quad + 12x_1^2 \sum_{\mathrm{a},\mathrm{b}\in D_k} p_{\mathrm{a}} p_{\mathrm{b}} (a_1-b_1)^2 (a_2-b_2)^2
\end{aligned}$$

$$= C_5(k) + C_6(k)x_1^2,$$

where, as above, $C_5(k)$ and $C_6(k)$ are both polynomials in $m_i(k)$ and $m_{i,j}(k)$, each of which converges to the corresponding polynomial in $m_i$ and $m_{i,j}$ as $k \to \infty$. Similar calculations show that $T_{1,3}(3) := 3\sum_{\mathrm{a,b} \in D_k} p_{\mathrm{a}} p_{\mathrm{b}} (a_2 - b_2)^2[(a_2 - b_2)^2 + 4x_2(a_2 - b_2) + 4x_2^2][2x_1(a_1 - b_1) + (a_1 - b_1)^2] = C_5(k) + C_6(k)x_2^2$.

Finally, calculations similar to (11) yield $T_2(3) := \sum_{\mathrm{a} \notin D_k \text{ or } \mathrm{b} \notin D_k} (\|\mathrm{x} + \mathrm{a} - \mathrm{b}\|_2^2 - \|\mathrm{x}\|_2^2)^3 \mathbb{P}\{X_{\mathrm{x}}^1 = \mathrm{x} + \mathrm{a}, X_{\mathrm{o}}^1 = \mathrm{b} \mid X_{\mathrm{o}}^0 = 0, X_{\mathrm{x}}^1 = \mathrm{x}\} = o(k^{-2})$ as $k \to \infty$. This establishes (5) and completes the proof of Lemma 3.1. □

3.3. $d \geq 4$. For notational simplicity we present the proof only for $d = 4$. Throughout this section the letters $\mathbf{u}$, $\mathbf{v}$ in bold font denote vectors in $\mathbb{Z}^4$, u, v in roman font denote vectors in $\mathbb{Z}^3$ and $u$, $v$ in italic font denote integers. We first show that on $\mathbb{Z}^4$, the graph $\mathcal{G}$ admits two distinct trees with positive probability, that is,

$$\mathbb{P}\{\mathcal{G} \text{ is disconnected}\} > 0. \tag{12}$$

Consider a random vector $X \in \mathbb{Z}^3$ defined as follows: for $k \geq 0$, let $\Delta_k := \{\mathrm{v} \in \mathbb{Z}^3 : \|\mathrm{v}\|_1 \leq k\}$ denote the three-dimensional diamond of radius $k$ and let $\delta\Delta_k := \{\mathrm{v} \in \mathbb{Z}^3 : \|\mathrm{v}\|_1 = k\}$ denote its boundary. As in (6), the distribution of the random vector $X$ is given by

$$\mathbb{P}(X = \mathrm{v}) = \begin{cases} p, & \text{if } \mathrm{v} = \mathrm{o}, \\ \dfrac{(1-p)^{\#\Delta_{k-1}}(1 - (1-p)^{\#\delta\Delta_k})}{\#\delta\Delta_k}, & \text{for } \mathrm{v} \in \delta\Delta_k,\ k \geq 1, \end{cases} \tag{13}$$

where $\mathrm{o} := (0,0,0)$ and $\#A$ denotes the cardinality of the set $A$. It may easily be checked that $\sum_{\mathrm{v} \in \mathbb{Z}^3} P(X = \mathrm{v}) = 1$.

Next, for a fixed vector $\mathbf{u} := (u(1), \dots, u(4)) \in \mathbb{Z}^4$, consider the graph $\mathcal{H} := (\mathcal{V} \cup \{\mathbf{u}\}, \mathcal{E} \cup \{\langle \mathbf{u}, h(\mathbf{u})\rangle\})$. For $n \geq 0$, let $h^n(\mathbf{u}) := (g^n(\mathbf{u}), t)$ for $g^n(\mathbf{u}) \in \mathbb{Z}^3$ and $t = u(4) - n \in \mathbb{Z}$. Here we take $h^0(\mathbf{u}) = \mathbf{u}$. Observe that for fixed $\mathbf{u}$, $g^n(\mathbf{u})$ has the same distribution as $(u(1), u(2), u(3)) + \sum_{i=1}^n X_i$, where $X_1, X_2, \dots$ are i.i.d. copies of $X$. Hence $\{g^n(\mathbf{u}) : n \geq 0\}$ is a symmetric random walk starting at $g^0(\mathbf{u}) = (u(1), u(2), u(3))$, with i.i.d. steps, each step size having distribution $X$. However, for $\mathbf{v} \in \mathbb{Z}^4$ with $v(4) = u(4)$, in the graph $(\mathcal{V} \cup \{\mathbf{u}, \mathbf{v}\}, \mathcal{E} \cup \{\langle \mathbf{u}, h(\mathbf{u})\rangle, \langle \mathbf{v}, h(\mathbf{v})\rangle\})$ the processes $\{g^n(\mathbf{u})\}_{n \geq 0}$ and $\{g^n(\mathbf{v})\}_{n \geq 0}$ are not independent and so, to obtain our theorem, we cannot use the fact that, with positive probability, two independent random walks on $\mathbb{Z}^3$ do not intersect. Nonetheless, if $\mathbf{u}$ and $\mathbf{v}$ are sufficiently far apart, their dependence on each other is weak. In the remainder of this section we formalize this notion of weak dependence by coupling two independent random walks and the processes $\{g^n(\mathbf{u}), g^n(\mathbf{v}) : n \geq 0\}$ and obtain the desired result.

For $\mathbf{v} = (\mathrm{v}, 0)$, given $\varepsilon > 0$ define the event

$$A_{n,\varepsilon}(\mathbf{v}) := \{g^{n^4}(\mathbf{v}) \in g^{n^4}(\mathbf{0}) + (\Delta_{n^{2(1+\varepsilon)}} \setminus \Delta_{n^{2(1-\varepsilon)}}), \quad g^i(\mathbf{v}) \neq g^i(\mathbf{0}) \text{ for all } i = 1, \dots, n^4\}, \tag{14}$$

where $\mathbf{0} := (0,0,0,0)$.

LEMMA 3.2. *For $0 < \varepsilon < 1/3$, there exist constants $C, \beta > 0$ and $n_0 \geq 1$ such that, for all $n \geq n_0$,*

$$\inf_{g^0(\mathbf{v}) \in \Delta_{n^{1+\varepsilon}} \setminus \Delta_{n^{1-\varepsilon}}} \mathbb{P}(A_{n,\varepsilon}(\mathbf{v})) \geq 1 - Cn^{-\beta}.$$

Assuming the above lemma, we proceed to complete the proof of (12). We shall return to the proof of the lemma later.

For $i \geq 1$ and $n \geq n_0$, let $\tau_i (= \tau_i(n)) := 1 + n^4 + (n^4)^2 + \cdots + (n^4)^{2^{i-1}}$ and take $\tau_0 = 1$. For fixed $\mathbf{v}$, we define

$$B_0 = B_0(\mathbf{v}) := \{g(\mathbf{v}) \in g(\mathbf{0}) + (\Delta_{n^{1+\varepsilon}} \setminus \Delta_{n^{1-\varepsilon}})\},$$

and having defined $B_0, \dots, B_{i-1}$, we define

$$B_i = B_i(\mathbf{v}) := \{g^{\tau_i}(\mathbf{v}) \in g^{\tau_i}(\mathbf{0}) + (\Delta_{n^{2^i(1+\varepsilon)}} \setminus \Delta_{n^{2^i(1-\varepsilon)}}) \text{ and } g^j(\mathbf{v}) \neq g^j(\mathbf{0}) \text{ for all } \tau_{i-1} + 1 \leq j \leq \tau_i\}.$$

Clearly,

$$\begin{aligned} &\mathbb{P}\{g^j(\mathbf{v}) \neq g^j(\mathbf{0}) \text{ for all } j \geq 1\} \\ &\quad \geq \mathbb{P}\left(\bigcap_{i=0}^{\infty} B_i\right) \\ &\quad = \lim_{i \to \infty} \mathbb{P}\left(\bigcap_{j=0}^{i} B_j\right) \\ &\quad = \lim_{i \to \infty} \prod_{l=1}^{i} \mathbb{P}\left(B_l \,\middle|\, \bigcap_{j=0}^{l-1} B_j\right) \mathbb{P}(B_0). \end{aligned} \tag{15}$$

Since $\mathbb{P}(B_0) > 0$, from (15) we have that $\mathbb{P}(g^j(\mathbf{v}) \neq g^j(\mathbf{0}) \text{ for all } j \geq 1) > 0$ if $\sum_{l=1}^{\infty} 1 - \mathbb{P}(B_l | \bigcap_{j=0}^{l-1} B_j) < \infty$.

For fixed $l \geq 1$, let $\mathbf{u}_1 := h^{\tau_l}(\mathbf{0})$ and $\mathbf{v}_1 := h^{\tau_l}(\mathbf{v})$. Now $\{(h^n(\mathbf{0}), h^n(\mathbf{v})) : n \geq 0\}$ being a Markov process and, since $g^0(\mathbf{v}_1)(\omega, \xi) \in g^0(\mathbf{u}_1)(\omega, \xi) + (\Delta_{n^{2^l(1+\varepsilon)}} \setminus \Delta_{n^{2^l(1-\varepsilon)}})$ for $(\omega, \xi) \in B_l(\mathbf{v})$, we have

$$\mathbb{P}\left(B_{l+1} \,\middle|\, \bigcap_{j=0}^{l} B_j\right)$$

$$
\begin{aligned}
&\geq \inf_1 \mathbb{P}\Big\{g^{(n^4)^{2^l}}(\mathbf{v}_1) \in g^{(n^4)^{2^l}}(\mathbf{u}_1) + (\Delta_{n^{2^{l+1}}(1+\varepsilon)} \setminus \Delta_{n^{2^{l+1}}(1-\varepsilon)}),\\
&\qquad\qquad g^k(\mathbf{u}_1) \neq g^k(\mathbf{v}_1) \text{ for all } k = 1, 2, \dots, (n^4)^{2^l}\Big\} \qquad (16)\\
&= \inf_2 \mathbb{P}(A_{n^{2^l},\varepsilon}(\mathbf{u})) \ \geq 1 - C(n^{2^l})^{-\beta},
\end{aligned}
$$

where $\inf_1$ is the infimum over all $\mathbf{u}_1, \mathbf{v}_1 \in \mathbb{Z}^4$ with $g^0(\mathbf{v}_1) \in g^0(\mathbf{u}_1) + (\Delta_{n^{2^l}(1+\varepsilon)} \setminus \Delta_{n^{2^l}(1-\varepsilon)})$ and $\inf_2$ is the infimum over all $\mathbf{u}$ with $g^0(\mathbf{u}) \in (\Delta_{n^{2^l}(1+\varepsilon)} \setminus \Delta_{n^{2^l}(1-\varepsilon)})$ and the last inequality follows from Lemma 3.2. Thus $\sum_{l=1}^{\infty}(1-\mathbb{P}(B_l | \bigcap_{j=0}^{l-1} B_j)) \leq C\sum_{l=1}^{\infty}(n^{2^l})^{-\beta} < \infty$, thereby completing the proof of (12).

To prove Lemma 3.2, we have to compare the trees $\{h^n(\mathbf{0})\}$ and $\{h^n(\mathbf{v})\}$ and independent "random walks" $\{\mathbf{0}+(\sum_{i=1}^n X_i, -n)\}$ and $\{\mathbf{v}+(\sum_{i=1}^n Y_i, -n)\}$, where $\{X_1, X_2, \dots\}$ and $\{Y_1, Y_2, \dots\}$ are independent collections of i.i.d. copies of the random variable $X$ given in (13).

We now describe a method to couple the trees and the independent random walks. Before embarking on the formal details of the coupling procedure, we present the main idea.

From a vertex $\mathbf{0}$ we construct the "path" $\{\mathbf{0} + (\sum_{i=1}^n X_i, -n)\}$. Now consider the vertex $\mathbf{v}$ with $\mathbf{v} = (v_1, v_2, v_3, 0)$. In case the diamond $D := \{\mathrm{u} \in \mathbb{Z}^3 : \|\mathrm{u}\|_1 \leq \|X_1\|_1\}$ is disjoint from the diamond $D' := \{\mathrm{u} \in \mathbb{Z}^3 : \|\mathrm{u} - (v_1, v_2, v_3)\|_1 \leq \|Y_1\|_1\}$, then we take $h^1(\mathbf{v}) = \{\mathbf{v} + (Y_1, -1)\}$. If the two diamonds are not disjoint, then we have to define $h^1(\mathbf{v})$ taking into account the configuration inside the diamond $D$. Similarly, we may obtain $h^2(\mathbf{v})$ by considering the diamonds $\{\mathrm{u} \in \mathbb{Z}^3 : \|\mathrm{u} - X_1\|_1 \leq \|X_2\|_1\}$ and $\{\mathrm{u} \in \mathbb{Z}^3 : \|\mathrm{u} - g^1(\mathbf{v})\|_1 \leq \|Y_2\|_1\}$. Note that if, for each $i = 1, \dots, n$, the two diamonds involved at the $i$th stage are disjoint, then the growth of the tree $\{(h^i(\mathbf{0}), h^i(\mathbf{v})) : 0 \leq i \leq n\}$ is stochastically equivalent to that of the pair of independent "random walks" $(\mathbf{0} + (\sum_{i=1}^n X_i, -n), \mathbf{v} + (\sum_{i=1}^n Y_i, -n))$.

We start with two vertices $\mathbf{u} := (\mathrm{u}, 0)$ and $\mathbf{v} := (\mathrm{v}, 0)$ in $\mathbb{Z}^4$ with $\mathrm{u}, \mathrm{v} \in \mathbb{Z}^3$. Let $\{U_1^{\mathrm{u}}(\mathrm{z}) : \mathrm{z} \in \mathbb{Z}^3\}, \{U_2^{\mathrm{u}}(\mathrm{z}) : \mathrm{z} \in \mathbb{Z}^3\}$ and $\{U_1^{\mathrm{v}}(\mathrm{z}) : \mathrm{z} \in \mathbb{Z}^3\}, \{U_2^{\mathrm{v}}(\mathrm{z}) : \mathrm{z} \in \mathbb{Z}^3\}$ be four independent collections of i.i.d. random variables, each of these random variables being uniformly distributed on $[0, 1]$.

Let $k_{\mathrm{u}}$ and $l_{\mathrm{v}}$ be defined as

$$
\begin{aligned}
k_{\mathrm{u}} &:= \min\{k : U_1^{\mathrm{u}}(\mathrm{z}) < p \text{ for some } \mathrm{z} \in (\mathrm{u} + \Delta_k)\},\\
l_{\mathrm{v}} &:= \min\{l : U_1^{\mathrm{v}}(\mathrm{z}) < p \text{ for some } \mathrm{z} \in (\mathrm{v} + \Delta_l)\}.
\end{aligned}
$$

Now define $m_{\mathrm{v}}$ as

$$
\begin{aligned}
m_{\mathrm{v}} := \min\{m :\ &\text{either } U_1^{\mathrm{v}}(\mathrm{z}) < p \text{ for some } \mathrm{z} \in (\mathrm{v} + \Delta_m) \setminus (\mathrm{u} + \Delta_{k_{\mathrm{u}}})\\
&\text{or } U_1^{\mathrm{u}}(\mathrm{z}) < p \text{ for some } \mathrm{z} \in (\mathrm{v} + \Delta_m) \cap (\mathrm{u} + \Delta_{k_{\mathrm{u}}})\}.
\end{aligned}
$$

Also, define the sets

$$\begin{aligned}N_{\mathrm{u}} &:= \{\mathrm{z} \in (\mathrm{u} + \Delta_{k_{\mathrm{u}}}) : U_1^{\mathrm{u}}(\mathrm{z}) < p\},\\ N_{\mathrm{v}}^1 &:= \{\mathrm{z} \in (\mathrm{v} + \Delta_{l_{\mathrm{v}}}) : U_1^{\mathrm{v}}(\mathrm{z}) < p\},\\ N_{\mathrm{v}}^2 &:= \{\mathrm{z} \in (\mathrm{v} + \Delta_{m_{\mathrm{v}}}) \setminus (\mathrm{u} + \Delta_{k_{\mathrm{u}}}) : U_1^{\mathrm{v}}(\mathrm{z}) < p\}\\ &\quad \cup \{\mathrm{z} \in (\mathrm{v} + \Delta_{m_{\mathrm{v}}}) \cap (\mathrm{u} + \Delta_{k_{\mathrm{u}}}) : U_1^{\mathrm{u}}(\mathrm{z}) < p\}.\end{aligned}$$

We pick:

(a) $\phi(\mathrm{u}) \in N_{\mathrm{u}}$ such that $U_2^{\mathrm{u}}(\phi(\mathrm{u})) = \min\{U_2^{\mathrm{u}}(\mathrm{z}) : \mathrm{z} \in N_{\mathrm{u}}\}$;
(b) $\zeta(\mathrm{v}) \in N_{\mathrm{v}}^1$ such that $U_2^{\mathrm{v}}(\zeta(\mathrm{v})) = \min\{U_2^{\mathrm{v}}(\mathrm{z}) : \mathrm{z} \in N_{\mathrm{v}}^1\}$;
(c) $\psi(\mathrm{v}) \in N_{\mathrm{v}}^2$ such that $U_2^{\mathrm{v}}(\psi(\mathrm{v})) = \min\{U_2^{\mathrm{v}}(\mathrm{z}) : \mathrm{z} \in N_{\mathrm{v}}^2\}$.

Taking $\phi^0(\mathrm{u}) = \mathrm{u}$, $\phi^n(\mathrm{u}) = \phi(\phi^{n-1}(\mathrm{u}))$, and similarly for $\zeta^n(\mathrm{v})$ and $\psi^n(\mathrm{v})$, we note that the distribution of $\{((\phi^n(\mathrm{u}), -n), (\zeta^n(\mathrm{v}), -n)) : n \geq 0\}$ is the same as that of $\{((\mathrm{u} + \sum_{i=1}^n X_i, -n), (\mathrm{v} + \sum_{i=1}^n Y_i, -n)) : n \geq 0\}$, that is, two independent "random walks," one starting from $(\mathrm{u}, 0)$ and the other starting from $(\mathrm{v}, 0)$. Also the distribution of $\{(h_n(\mathrm{u}, 0), h_n(\mathrm{v}, 0)) : n \geq 0\}$ and that of $\{((\phi^n(\mathrm{u}), -n), (\psi^n(\mathrm{v}), -n)) : n \geq 0\}$ are identical. Thus, the procedure described above may be used to construct the trees from $(\mathrm{u}, 0)$ and $(\mathrm{v}, 0)$.

Now observe that $\{(\phi^n(\mathrm{u}), -n)\}$ describes both the random walk and the tree starting from $(\mathrm{u}, 0)$. Also if $\Delta_{k_{\mathrm{u}}} \cap \Delta_{m_{\mathrm{v}}} = \varnothing$, then $m_{\mathrm{v}} = l_{\mathrm{v}}$ and, more importantly, $\zeta(\mathrm{v}) = \psi(\mathrm{v})$. Hence the "random walk" and the tree from $(\mathrm{u}, 0)$ are coupled and so are the "random walk" and the tree from $(\mathrm{v}, 0)$. In particular, this happens when both $k_{\mathrm{u}} < [\|\mathrm{u} - \mathrm{v}\|_1/2]$ and $m_{\mathrm{v}} < [\|\mathrm{u} - \mathrm{v}\|_1/2]$. Let $k_0 = \|\mathrm{u} - \mathrm{v}\|_1/2$. From the above discussion, we have

$$\begin{aligned}\mathbb{P}(\{\zeta(\mathrm{v}) \neq \psi(\mathrm{v})\}) &\leq \mathbb{P}(\{(U_1^{\mathrm{u}}(\mathrm{z})) > p \text{ for all } \mathrm{z} \in (\mathrm{u} + \Delta_{k_0})\}\\ &\qquad \cup \{(U_1^{\mathrm{v}}(\mathrm{z})) > p \text{ for all } \mathrm{z} \in (\mathrm{v} + \Delta_{k_0})\})\\ &= 2\mathbb{P}(\{(U_1^{\mathrm{u}}(\mathrm{z})) > p \text{ for all } \mathrm{z} \in (\mathrm{u} + \Delta_{k_0})\})\\ &= 2(1-p)^{\#\Delta_{k_0}}.\end{aligned}$$

Since $(1/2)k^3 \leq \#\Delta_k \leq 2k^3$, the above inequality gives

$$\mathbb{P}(\{\zeta(\mathrm{v}) = \psi(\mathrm{v})\}) \geq 1 - C_1 \exp(-C_2 \|\mathrm{u} - \mathrm{v}\|_1^3) \tag{17}$$

for constants $C_1 = 2$ and $C_2 = (1/2)|\log(1-p)|$.

With the above estimate at hand, we look at the process $\{(\phi^n(\mathrm{u}), \zeta^n(\mathrm{v})) : n \geq 0\}$. Without loss of generality we take $\mathrm{u} = \mathrm{o}$. For $\varepsilon > 0$ and constant $K > 0$ (to be specified later), define

$$\begin{aligned}B_{n,\varepsilon}(\mathrm{v}) := \{&\zeta^{n^4}(\mathrm{v}) \in \phi^{n^4}(\mathrm{o}) + (\Delta_{n^{2(1+\varepsilon)}} \setminus \Delta_{n^{2(1-\varepsilon)}}),\\ &\|\zeta^i(\mathrm{v}) - \phi^i(\mathrm{o})\|_1 \geq K \log n \text{ for all } i = 1, \ldots, n^4\}.\end{aligned} \tag{18}$$

This event is an independent random walk version of the event $A_{n,\varepsilon}(\mathrm{v},0)$ defined in (14), except that here we require that the two random walks come no closer than $K\log n$ at any stage.

We will show that there exists $\alpha>0$ such that

$$\sup_{\mathrm{v}\in(\Delta_{n^{(1+\varepsilon)}}\setminus\Delta_{n^{(1-\varepsilon)}})}\mathbb{P}((B_{n,\varepsilon}(\mathrm{v}))^c)<C_3n^{-\alpha}\tag{19}$$

for some constant $C_3>0$.

Since $(B_{n,\varepsilon}(\mathrm{v}))^c\subseteq E_{n,\varepsilon}(\mathrm{v})\cup F_{n,\varepsilon}(\mathrm{v})\cup G_{n,\varepsilon}(\mathrm{v})$, where

$$E_{n,\varepsilon}(\mathrm{v}):=\{\|\zeta^i(\mathrm{v})-\phi^i(\mathrm{o})\|_1\le K\log n \text{ for some } i=1,\dots,n^4\},$$

$$F_{n,\varepsilon}(\mathrm{v}):=\{\zeta^{n^4}(\mathrm{v})\notin\phi^{n^4}(\mathrm{o})+\Delta_{n^{2(1+\varepsilon)}}\},$$

$$G_{n,\varepsilon}(\mathrm{v}):=\{\zeta^{n^4}(\mathrm{v})\in\phi^{n^4}(\mathrm{o})+\Delta_{n^{2(1-\varepsilon)}}\},$$

to prove (19) it suffices to show the following.

LEMMA 3.3. *There exist $\alpha>0$ and constants $C_4,C_5,C_6>0$ such that, for all $n$ sufficiently large, we have:*

(a) $\sup_{\mathrm{v}\in(\Delta_{n^{(1+\varepsilon)}}\setminus\Delta_{n^{(1-\varepsilon)}})}\mathbb{P}(E_{n,\varepsilon}(\mathrm{v}))<C_4n^{-\alpha}$,
(b) $\sup_{\mathrm{v}\in(\Delta_{n^{(1+\varepsilon)}}\setminus\Delta_{n^{(1-\varepsilon)}})}\mathbb{P}(F_{n,\varepsilon}(\mathrm{v}))<C_5n^{-\alpha}$,
(c) $\sup_{\mathrm{v}\in(\Delta_{n^{(1+\varepsilon)}}\setminus\Delta_{n^{(1-\varepsilon)}})}\mathbb{P}(G_{n,\varepsilon}(\mathrm{v}))<C_6n^{-\alpha}$.

PROOF. First we fix $\mathrm{v}\in(\Delta_{n^{(1+\varepsilon)}}\setminus\Delta_{n^{(1-\varepsilon)}})$. Since $\{(\phi^n(\mathrm{o}),\zeta^n(\mathrm{v})):n\ge0\}$ and $\{(\sum_{i=1}^n X_i,\mathrm{v}+\sum_{i=1}^n Y_i):n\ge0\}$ have the same distribution, we have

$$\begin{aligned}\mathbb{P}(E_{n,\varepsilon}(\mathrm{v}))&=\mathbb{P}\left\{\left\|\sum_{j=1}^i X_j-\left(\mathrm{v}+\sum_{j=1}^i Y_j\right)\right\|_1\le K\log n \text{ for some } i=1,\dots,n^4\right\}\\&=\mathbb{P}\left\{\sum_{j=1}^i X_j-\sum_{j=1}^i Y_j\in(\mathrm{v}+\Delta_{K\log n}) \text{ for some } i=1,\dots,n^4\right\}\\&\le\mathbb{P}\left\{\sum_{j=1}^i X_j-\sum_{j=1}^i Y_j\in(\mathrm{v}+\Delta_{K\log n}) \text{ for some } i\ge1\right\}\\&=\mathbb{P}\left(\bigcup_{\mathrm{z}\in(\mathrm{v}+\Delta_{K\log n})}\left\{\sum_{j=1}^i X_j-\sum_{j=1}^i Y_j=\mathrm{z} \text{ for some } i\ge1\right\}\right).\end{aligned}$$

Now $\sum_{j=1}^i(X_j-Y_j)$ is an aperiodic, isotropic, symmetric random walk whose steps are i.i.d. with each step having the same distribution as $X-Y$, where $Y$ is an independent copy of $X$. Since $\mathrm{Var}(X-Y)=2\,\mathrm{Var}(X)=2\sigma^2I$ [where $\sigma^2:=\mathrm{Var}(X(1))$ and $\mathrm{Var}(X)$ denotes the variance-covariance matrix of $X$] and $\sum_{\mathrm{u}\in\mathbb{Z}^3}|\mathrm{u}|^2\mathbb{P}(X-Y=\mathrm{u})<\infty$, by Proposition P26.1 of Spitzer [(1964), page 308],

$$\lim_{|\mathrm{z}|\to\infty}|\mathrm{z}|\,\mathbb{P}\left\{\sum_{j=1}^{i}X_j-\sum_{j=1}^{i}Y_j=\mathrm{z}\text{ for some }i\ge 1\right\}=(4\pi\operatorname{Var}(X(1)))^{-1}. \tag{20}$$

For $\mathrm{v}\in(\Delta_{n^{(1+\varepsilon)}}\setminus\Delta_{n^{(1-\varepsilon)}})$ and $\mathrm{z}\in\mathrm{v}+\Delta_{K\log n}$, we must have that, for all $n$ sufficiently large, $|\mathrm{z}|\ge n^{1-\varepsilon}/2$. Thus for all $n$ sufficiently large and for some constants $C_7,C_8,C_9>0$, we have, using (20),

$$\begin{aligned}\mathbb{P}(E_{n,\varepsilon}(\mathrm{v}))&\le\sum_{\mathrm{z}\in(\mathrm{v}+\Delta_{K\log n})}\mathbb{P}\left\{\sum_{j=1}^{i}X_j-\sum_{j=1}^{i}Y_j=\mathrm{z}\text{ for some }i\ge 1\right\}\\&\le C_7(K\log n)^3C_8(n^{-(1-\varepsilon)})\\&\le C_9n^{-(1-\varepsilon/2)}.\end{aligned}$$

This completes the proof of Lemma 3.3(a).

For the next part of the lemma, observe that, for sufficiently large $n$ and all $\mathrm{v}\in\Delta_{n^{(1+\varepsilon)}}\setminus\Delta_{n^{(1-\varepsilon)}}$,

$$\begin{aligned}\mathbb{P}(F_{n,\varepsilon}(\mathrm{v}))&=\mathbb{P}\left\{\mathrm{v}+\sum_{j=1}^{n^4}(X_j-Y_j)\notin\Delta_{n^{2(1+\varepsilon)}}\right\}\\&=\mathbb{P}\left\{\left\|\mathrm{v}+\sum_{j=1}^{n^4}(X_j-Y_j)\right\|_1>n^{2(1+\varepsilon)}\right\}\\&\le\mathbb{P}\left\{\left\|\sum_{j=1}^{n^4}(X_j-Y_j)\right\|_1>n^{2(1+\varepsilon)}-n^{(1+\varepsilon)}\right\}\\&\le\mathbb{P}\left\{\left\|\sum_{j=1}^{n^4}(X_j-Y_j)\right\|_1>n^{2(1+\varepsilon)}/2\right\}.\end{aligned}\tag{21}$$

To estimate the above probability let $X-Y=Z=(Z(1),Z(2),Z(3))$, where $E(Z(i))=0$ and $\operatorname{Var}(Z(i))=2\sigma^2$. Then, letting $\sum_{j=1}^{k}(X_j-Y_j)(i)$ denote the $i$th co-ordinate of the process $\sum_{j=1}^{k}(X_j-Y_j)$ and using Chebyshev's inequality, we have

$$\begin{aligned}&\mathbb{P}\left\{\left\|\sum_{j=1}^{n^4}(X_j-Y_j)\right\|_1>\frac{n^{2(1+\varepsilon)}}{2}\right\}\\&\qquad\le\mathbb{P}\left\{\bigcup_{i=1}^{3}\left\{\left|\sum_{j=1}^{n^4}(X_j-Y_j)(i)\right|>\frac{n^{2(1+\varepsilon)}}{6}\right\}\right\}\end{aligned}$$

$$\leq 3\mathbb{P}\left\{\left|\sum_{j=1}^{n^4}(X_j - Y_j)(1)\right| > \frac{n^{2(1+\varepsilon)}}{6}\right\}$$
$$\leq \frac{3n^4 \operatorname{Var}(Z(1))}{(n^{2(1+\varepsilon)}/6)^2}$$
$$\leq \frac{C_{11}}{n^{4\varepsilon}},$$

for some constant $C_{11} > 0$. Combining the above inequality with that obtained in (21), we have

$$\sup_{\mathrm{v}\in(\Delta_{n^{(1+\varepsilon)}}\setminus\Delta_{n^{(1-\varepsilon)}})} \mathbb{P}(F_{n,\varepsilon}(\mathrm{v})) \leq \frac{C_{11}}{n^{4\varepsilon}},$$

which proves Lemma 3.3(b).

Finally, for the last part of the lemma, we have that if $0 < \varepsilon < 1/3$ and $\mathrm{v} \in \Delta_{n^{(1+\varepsilon)}} \setminus \Delta_{n^{(1-\varepsilon)}}$, for all sufficiently large $n$, $\|\mathrm{v}\|_1 < n^{2(1-\varepsilon)}$. Therefore,

$$\begin{aligned}
\mathbb{P}(G_{n,\varepsilon}(\mathrm{v})) &\leq \mathbb{P}\left\{\left\|\mathrm{v} + \sum_{j=1}^{n^4}(X_j - Y_j)\right\|_1 < n^{2(1-\varepsilon)}\right\} \\
&\leq \mathbb{P}\left\{\left\|\sum_{j=1}^{n^4}(X_j - Y_j)\right\|_1 < \|\mathrm{v}\|_1 + n^{2(1-\varepsilon)}\right\} \\
&\leq \mathbb{P}\left\{\left\|\sum_{j=1}^{n^4}(X_j - Y_j)\right\|_1 < 2n^{2(1-\varepsilon)}\right\} \qquad (22)\\
&\leq \mathbb{P}\left\{\bigcup_{i=1}^{3}\left\{\left|\sum_{j=1}^{n^4}(X_j - Y_j)(i)\right| < \frac{2n^{2(1-\varepsilon)}}{3}\right\}\right\} \\
&\leq 3\mathbb{P}\left\{\frac{|\sum_{j=1}^{n^4}(X_j - Y_j)(1)|}{n^2} < \frac{2n^{-2\varepsilon}}{3}\right\}.
\end{aligned}$$

By the central limit theorem, as $n \to \infty$, $\sum_{j=1}^{n^4}(X_j - Y_j)(1)/(\sqrt{2}\sigma n^2)$ converges in distribution to a random variable $N$ (say) with a standard normal distribution. Thus

$$\begin{aligned}
&\mathbb{P}\left\{\frac{|\sum_{j=1}^{n^4}(X_j - Y_j)(1)|}{n^2} < \frac{2n^{-2\varepsilon}}{3}\right\} - \mathbb{P}\left\{|N| < \frac{\sqrt{2}n^{-2\varepsilon}}{3\sigma}\right\} \qquad (23)\\
&\quad\leq \left|\mathbb{P}\left\{\frac{|\sum_{j=1}^{n^4}(X_j - Y_j)(1)|}{\sqrt{2}\sigma n^2} < \frac{\sqrt{2}n^{-2\varepsilon}}{3\sigma}\right\} - \mathbb{P}\left\{|N| < \frac{\sqrt{2}n^{-2\varepsilon}}{3\sigma}\right\}\right|.
\end{aligned}$$

Of the terms in the above inequality, we have

$$\begin{aligned}\mathbb{P}\left\{|N| \le \frac{\sqrt{2}n^{-2\varepsilon}}{3\sigma}\right\} &= \int_{-\sqrt{2}n^{-2\varepsilon}(3\sigma)^{-1}}^{\sqrt{2}n^{-2\varepsilon}(3\sigma)^{-1}} \frac{1}{\sqrt{2\pi}} \exp\left(-\frac{x^2}{2}\right) dx \\ &\le \frac{2\sqrt{2}n^{-2\varepsilon}(3\sigma)^{-1}}{\sqrt{2\pi}},\end{aligned} \tag{24}$$

and we use Berry–Essen bounds [see Chow and Teicher (1978), Corollary 9.4, page 300], to obtain

$$\begin{aligned}&\left|\mathbb{P}\left\{\frac{|\sum_{j=1}^{n^4}(X_j - Y_j)(1)|}{\sqrt{2}\sigma n^2} < \frac{\sqrt{2}n^{-2\varepsilon}}{3\sigma}\right\} - \mathbb{P}\left\{|N| < \frac{\sqrt{2}n^{-2\varepsilon}}{3\sigma}\right\}\right| \\ &\quad \le \left|\mathbb{P}\left\{\frac{\sum_{j=1}^{n^4}(X_j - Y_j)(1)}{\sqrt{2}\sigma n^2} < \frac{\sqrt{2}n^{-2\varepsilon}}{3\sigma}\right\} - \mathbb{P}\left\{N < \frac{\sqrt{2}n^{-2\varepsilon}}{3\sigma}\right\}\right| \\ &\qquad + \left|\mathbb{P}\left\{\frac{\sum_{j=1}^{n^4}(X_j - Y_j)(1)}{\sqrt{2}\sigma n^2} \le -\frac{\sqrt{2}n^{-2\varepsilon}}{3\sigma}\right\} - \mathbb{P}\left\{N \le -\frac{\sqrt{2}n^{-2\varepsilon}}{3\sigma}\right\}\right| \\ &\quad \le 2\sup_{x\in\mathbb{R}}\left|\mathbb{P}\left\{\frac{\sum_{j=1}^{n^4}(X_j - Y_j)(1)}{\sqrt{2}\sigma n^2} \le x\right\} - P\{N \le x\}\right| \\ &\quad \le \frac{C_{12}E(Z_1^4)}{n^4\sigma^4},\end{aligned} \tag{25}$$

for some constant $C_{12} > 0$. Combining (22)–(24), we have Lemma 3.3(c). □

Proof of Lemma 3.2. Let $\mathbf{v} := (\mathrm{v}, 0) \in \mathbb{Z}^4$. Observe that $A_{n,\varepsilon}(\mathbf{v}) \supseteq B_{n,\varepsilon}(\mathrm{v}) \cap \{g^i(\mathbf{0}) = \sum_{j=1}^i X_j,\ g^i(\mathbf{v}) = \mathrm{v} + \sum_{j=1}^i Y_j$ for all $1 \le i \le n^4\}$. Hence

$$\begin{aligned}&\mathbb{P}(A_{n,\varepsilon}(\mathbf{v})) \\ &\quad \ge \mathbb{P}\left\{B_{n,\varepsilon}(\mathrm{v}) \cap \left\{g^i(\mathbf{0}) = \sum_{j=1}^i X_j,\ g^i(\mathbf{v}) = \mathrm{v} + \sum_{j=1}^i Y_j \text{ for } 1 \le i \le n^4\right\}\right\} \\ &\quad = \mathbb{P}\left\{B_{n,\varepsilon}(\mathrm{v}) \cap \left\{g^i(\mathbf{0}) = \sum_{j=1}^i X_j,\ g^i(\mathbf{v}) = \mathrm{v} + \sum_{j=1}^i Y_j \text{ for } 1 \le i \le n^4 - 1\right\}\right\} \\ &\qquad \times \mathbb{P}\left\{g^{n^4}(\mathbf{0}) = \sum_{j=1}^{n^4} X_j, g^{n^4}(\mathbf{v}) = \mathrm{v} + \sum_{j=1}^{n^4} Y_j \right| \\ &\qquad\qquad B_{n,\varepsilon}(\mathrm{v}) \cap \left\{g^i(\mathbf{0}) = \sum_{j=1}^i X_j, g^i(\mathbf{v}) = \mathrm{v} + \sum_{j=1}^i Y_j \text{ for } 1 \le i \le n^4 - 1\right\}\Bigg\}\end{aligned}$$

$$\geq \mathbb{P}\left\{B_{n,\varepsilon}(\mathrm{v}) \cap \left\{g^i(\mathbf{0}) = \sum_{j=1}^{i} X_j, g^i(\mathbf{v}) = \mathrm{v} + \sum_{j=1}^{i} Y_j \text{ for } 1 \leq i \leq n^4 - 1\right\}\right\}$$
$$\times (1 - C_1 \exp(-C_2(K \log n)^3)),$$

where the last inequality follows from (17) after noting that, given $B_{n,\varepsilon}(\mathrm{v})$, $g^i(\mathbf{0}) = \sum_{j=1}^{i} X_j$ and $g^i(\mathbf{v}) = \mathrm{v} + \sum_{j=1}^{i} Y_j$ hold for all $1 \leq i \leq n^4 - 1$, we have $\|g^{n^4-1}(\mathbf{0}) - g^{n^4-1}(\mathbf{v})\|_1 \geq K \log n$. The above argument may be used iteratively for $i = 1, \ldots, n^4 - 1$, and together with (19), we have

$$\begin{aligned}
\mathbb{P}(A_{n,\varepsilon}(\mathbf{v})) &\geq (1 - C_1 \exp(-C_2(K \log n)^3))^{n^4} \mathbb{P}(B_{n,\varepsilon}(\mathrm{v})) \\
&\geq (1 - C_1 n^4 \exp(-C_2 K^3 \log n))(1 - C_3 n^{-\alpha}) \\
&\geq (1 - C_1 n^4 n^{-C_2 K^3})(1 - C_3 n^{-\alpha}) \\
&= (1 - C_1 n^{-C_2 K^3 + 4})(1 - C_3 n^{-\alpha}).
\end{aligned}$$

Taking $K$ such that $C_2 K^3 > 4$ [i.e., $K^3 > 8|\log(1-p)|^{-1}$], we have

$$\begin{aligned}
\mathbb{P}(A_{n,\varepsilon}(\mathbf{v})) &\geq 1 - C_1 n^{-C_2 K^3 + 4} - C_3 n^{-\alpha} \\
&\geq 1 - C n^{-\beta},
\end{aligned}$$

for some constant $C > 0$ and $\beta := \min\{\alpha, C_2 K^3 - 4\} > 0$. This completes the proof of Lemma 3.2. □

Finally, to complete the theorem we need to show that $\mathcal{G}$ admits infinitely many trees almost surely. For $k \geq 2$, define $D^k(n,\varepsilon) := \{(\mathbf{u}_1, \mathbf{u}_2, \ldots, \mathbf{u}_k) : \mathbf{u}_i \in \mathbb{Z}^4$ such that $n^{1-\varepsilon} \leq \|g^0(\mathbf{u}_i) - g^0(\mathbf{u}_j)\|_1 \leq n^{1+\varepsilon}$ for all $i \neq j\}$. Define the event $A(n, \varepsilon, \mathbf{u}_1, \mathbf{u}_2, \ldots, \mathbf{u}_k) := \{n^{2(1-\varepsilon)} \leq \|g^{n^4}(\mathbf{u}_i) - g^{n^4}(\mathbf{u}_j)\|_1 \leq n^{2(1+\varepsilon)}$ and $g^t(\mathbf{u}_i) \neq g^t(\mathbf{u}_j)$ for all $t = 1, \ldots, n^4$ and for all $i \neq j\}$. Using Lemma 3.2, we can easily show, for $0 < \varepsilon < 1/3$ and for all large $n$,

$$\inf\{\mathbb{P}(A(n, \varepsilon, \mathbf{u}_1, \mathbf{u}_2, \ldots, \mathbf{u}_k) : (\mathbf{u}_1, \mathbf{u}_2, \ldots, \mathbf{u}_k) \in D^k(n,\varepsilon))\} \geq 1 - \frac{C_k}{n^\beta}, \tag{26}$$

where $C_k$ is a constant independent of $n$ (depending on $k$) and $\beta$ is as in Lemma 3.2. We may now imitate the method following the statement of Lemma 3.2 to obtain

$$\mathbb{P}\{g^t(\mathbf{u}_i) \neq g^t(\mathbf{u}_j) \text{ for all } t \geq 1 \text{ and for } 1 \leq i \neq j \leq k\} > 0.$$

Thus, by translation invariance and ergodicity, we have that, for all $k \geq 2$,

$$\mathbb{P}\{\mathcal{G} \text{ contains at least } k \text{ trees}\} = 1.$$

This shows that $\mathcal{G}$ contains infinitely many trees almost surely.

**4. Geometry of the graph $\mathcal{G}$.** We now prove Theorem 2.2 for $d = 2$; with minor modifications the same argument carries through for any dimensions. The idea behind this proof was suggested by the referee.

For $t \in \mathbb{Z}$, consider the set $N_t := \mathcal{G} \cap \{y = t\}$, the set of open vertices on the line $\{y = t\}$. For $x \in N_t$ and $n \geq 0$, let $B_t^n(x) := \{y \in N_{t+n} : h^n(y) = x\}$ be the set of the $n$th-order ancestors of the vertex $x \in N_t$. Now consider the set of vertices in $N_t$ which have $n$th-order ancestors, that is, $M_t^{(n)} := \{x \in N_t : B_t^n(x) \neq \varnothing\}$. Clearly, $M_t^{(n)} \subseteq M_t^{(m)}$ for $n > m$ and so $R_t := \lim_{n\to\infty} M_t^{(n)} = \bigcap_{n\geq 0} M_t^{(n)}$ is well defined. Moreover, this is the set of vertices in $N_t$ which have bi-infinite paths. We want to show that $\mathbb{P}(R_t = \varnothing) = 1$ for all $t \in \mathbb{Z}$. Since $\{R_t : t \in \mathbb{Z}\}$ is stationary, it suffices to show that $\mathbb{P}(R_0 = \varnothing) = 1$.

First note that by the translation invariance of the model, $\mathbb{P}(\{\#R_0 = 0\} \cup \{\#R_0 = \infty\}) = 1$. Now suppose $\mathbb{P}(\#R_0 = \infty) > 0$. A vertex $x \in R_t$ is called a *branching point* if $\#(B_t^1(x) \cap R_{t+1}) \geq 2$, that is, $x$ has at least two distinct infinite branches of ancestors. Note that this notion of "branching point" is similar to that of "encounter point" of Burton and Keane (1989). As in their proof of the uniqueness of the percolation cluster, our proof essentially uses the fact that it is impossible to embed a tree in a lattice.

We first show that

$$\mathbb{P}(\text{Origin is a branching point}) > 0. \tag{27}$$

Since $\mathbb{P}(\#R_0 = \infty) > 0$, we may fix two vertices $x = (x_1, 1)$ and $y = (y_1, 1)$ such that

$$\mathbb{P}(x, y \in (B_0^1(\mathbf{0}) \cap R_1)) > 0.$$

Thus the event $E_1 := \{B_1^n(x) \neq \varnothing, B_1^n(y) \neq \varnothing \text{ for all } n \geq 1\}$ has positive probability. Further, this event depends only on sites $\{u := (u_1, u_2) : u_2 \geq 1\}$. Now, consider the event $E_2 := \{(i, 0)$ is closed for all $i \neq 0$ with $-2\max\{|x_1| + 1, |y_1| + 1\} \leq i \leq 2\max\{|x_1| + 1, |y_1| + 1\}$ and $(0, 0)$ is open$\}$. Clearly $\mathbb{P}(E_2) > 0$. Since $E_1$ and $E_2$ depend on disjoint sets of vertices, we have

$$\mathbb{P}(\text{Origin is a branching point}) \geq \mathbb{P}(E_1 \cap E_2) = \mathbb{P}(E_1)\mathbb{P}(E_2) > 0.$$

Now, we define $r_0(n) := \#(R_0 \cap ([-n, n] \times \{0\}))$ and $r_1(n) := \#(R_1 \cap ([-n, n] \times \{1\}))$. We arrange the points of $R_0 \cap ([-n, n] \times \{0\})$ as $u_1, \dots, u_{r_0(n)}$, in an increasing order of the $x$ coordinates. By our construction of $\mathcal{G}$, neither $u_2$ nor $u_{r_0(n)-1}$ nor any of the vertices between them can be connected to a vertex on $N_1$ which lies outside $[-n, n] \times \{1\}$. Thus, each of the vertices $u_2, u_3, \dots, u_{r_0(n)-1}$ will have at least one ancestor in the set $R_1 \cap ([-n, n] \times \{1\})$. Moreover, each of the branching points in $u_2, \dots, u_{r_0(n)-1}$ has at least two distinct ancestors in the set $R_1 \cap ([-n, n] \times \{1\})$. Thus, if $r_0^{(2)}(n)$ is the number of branching points in $[-n, n] \times \{1\}$, we must have

$$r_1(n) - (r_0(n) - 2) \geq r_0^{(2)}(n) - 2. \tag{28}$$

But, by stationarity, we have $\mathbb{E}(r_1(n)) = \mathbb{E}(r_0(n))$ for all $n \geq 1$. Thus, for $n$ sufficiently large, from (27) we have

$$\begin{aligned} 0 &= \mathbb{E}(r_1(n) - r_0(n)) \geq \mathbb{E}r_0^{(2)}(n) - 4 \\ &= (2n+1)\mathbb{P}(\text{Origin is a branching point}) - 4 > 0. \end{aligned}$$

This contradiction establishes Theorem 2.2.

**5. Limit theorem.** We first prove Theorem 2.3(a). The proof of the next part of the theorem is similar and thus omitted. For simplicity in notation we shall prove the result for $d = 2$; however, our method is also valid for higher dimensions.

Fix $\nu \geq 0$. Let $B_n := [1,n] \times [1,n]$ be a box of width $n$ and, for $(i,j) \in B_n \cap \mathbb{Z}^2$, define random variables $Y_{i,j}$ as

$$Y_{i,j} := \begin{cases} 1, & \text{if the degree of the vertex } (i,j) \text{ in } B_n \cap \mathcal{V} \text{ is } \nu + 1, \\ 0, & \text{otherwise.} \end{cases}$$

Note for a vertex $(i,j)$, $Y_{i,j} = 1$ if and only if there are exactly $\nu$ edges "going up" from $(i,j)$ and one edge "going down" from it.

Let $Y_j^{(n)} := \sum_{i=1}^{n} (Y_{i,j} - \mathbb{E}(Y_{i,j}))$ and $S_n := \sum_{j=1}^{n} Y_j^{(n)}$. To prove Theorem 2.3 we need to show that the distribution of $S_n/n$ is asymptotically normal.

Towards this end, first observe that, for fixed $j$, $\{Y_{i,j}\}_{i \geq 1}$ is an $\alpha$-mixing sequence of random variables; that is, for all $m \geq 1$, $A \in \sigma(Y_{1,j}, Y_{2,j}, \dots, Y_{m,j})$ and $B \in \sigma(Y_{m+n,j}, Y_{m+n+1,j}, \dots)$, we have $|\mathbb{P}(A \cap B) - \mathbb{P}(A)\mathbb{P}(B)| \leq \alpha_n$, where $\alpha_n \to 0$ as $n \to \infty$. Indeed, given $A$ and $B$ as above, define

$$\begin{aligned} E := \Big\{ &\text{there exists an open vertex in each of the sets} \\ &\Big\{(i,j) : m + \frac{n}{4} \leq i \leq m + \frac{3n}{8}\Big\}, \Big\{(i,j+1) : m + \frac{3n}{8} \leq i \leq m + \frac{n}{2}\Big\}, \\ &\Big\{(i,j+1) : m + \frac{n}{2} \leq i \leq m + \frac{5n}{8}\Big\}, \Big\{(i,j) : m + \frac{3n}{8} \leq i \leq m + \frac{3n}{4}\Big\}\Big\}. \end{aligned}$$

Now $\mathbb{P}(E) = (1 - (1-p)^{n/8})^4 \to 1$ as $n \to \infty$. Also, given $E$, the event $A$ depends only on the configuration of the vertices $\{(i, j-1) : i \leq m + \frac{n}{4}\}$, $\{(i,j) : i \leq m\}$ and $\{(i, j+1) : i < m + \frac{n}{2}\}$, while the event $B$ depends on the vertices $\{(i, j-1) : i \geq m + \frac{3n}{4}\}$, $\{(i,j) : i \geq m + n\}$ and $\{(i, j+1) : i > m + \frac{n}{2}\}$. These sets of vertices being disjoint, given $E$, $A$ and $B$ are conditionally independent, a simple conditioning argument now yields that, for $n$ large enough,

$$|\mathbb{P}(A \cap B) - \mathbb{P}(A)\mathbb{P}(B)| \leq 5\mathbb{P}(E^c) \leq C_1 \exp(-C_2 n) \quad \text{for constants } C_1, C_2 > 0. \tag{29}$$

Also observe that, for fixed $i$, $\{Y_{i,j}\}_{j\geq 1}$ is a one-dependent sequence of random variables; that is, for fixed $i$, $Y_{i,j}$ is independent of $Y_{i,j'}$ for $j' \neq j-1, j, j+1$.

Now, for some $0 < \delta < 1$ to be chosen later and for $0 \leq k < r_n$, where $r_n := \lfloor \frac{n}{\lfloor n^\delta \rfloor + 1} \rfloor$, let

$$\begin{aligned} W_{k+1}^{(n)} &:= Y_{k\lfloor n^\delta \rfloor + k+1}^{(n)} + \cdots + Y_{(k+1)\lfloor n^\delta \rfloor + k}^{(n)}, \\ \eta_{k+1}^{(n)} &:= Y_{(k+1)\lfloor n^\delta \rfloor + k+1}^{(n)}, \\ E_n &:= Y_{r_n(\lfloor n^\delta \rfloor + 1)+1}^{(n)} + \cdots + Y_n^{(n)}. \end{aligned}$$

First we show that, for any $r \geq 1$, there exists a constant $C > 0$ such that

$$\mathbb{E}(Y_1^{(n)} + \cdots + Y_r^{(n)})^4 \leq C r^2 n^2. \tag{30}$$

Indeed note that, as in the proof of the first part of Theorem 27.5 of Billingsley (1979), we have $\mathbb{E}(Y_i^{(n)})^4 = \mathbb{E}(Y_1^{(n)})^4 \leq K n^2$ for some constant $K > 0$. Now

$$\mathbb{E}\left(\sum_{k=1}^{r} Y_k^{(n)}\right)^4 = \sum_{k,l,s,t=1}^{r} \mathbb{E}(Y_k^{(n)} Y_l^{(n)} Y_s^{(n)} Y_t^{(n)}), \tag{31}$$

and using the fact that $\{Y_k^{(n)}\}_{k\geq 1}$ is a one-dependent sequence of random variables, the Cauchy–Schwarz inequality and that $\mathbb{E}Y_1^{(n)} = 0$, we obtain after some elementary calculations

$$\mathbb{E}\left(\sum_{k=1}^{r} Y_k^{(n)}\right)^4 \leq 2r\mathbb{E}(Y_1^{(n)})^4 + r^2 \mathbb{E}(Y_1^{(n)})^4.$$

Here the term $2r\mathbb{E}(Y_1^{(n)})^4$ comes from the terms in the sum $\sum_{j,k,s,t=1}^{r} \mathbb{E}(Y_j^{(n)} Y_k^{(n)} \times Y_s^{(n)} Y_t^{(n)})$ when $j, k, s, t$ are close to each other so as to have dependence among all the four random variables making the product, while the term $r^2\mathbb{E}(Y_1^{(n)})^4$ comes from the terms of the sum when $j, k$ are close to each other, $s, t$ are close to each other, but there is independence between $(Y_j^{(n)}, Y_k^{(n)})$ and $(Y_s^{(n)}, Y_t^{(n)})$. This proves (30).

Now taking $r = \lfloor n^\delta \rfloor$, and using the fact that $W_1^{(n)}, W_2^{(n)}, \ldots$ are i.i.d. random variables, we have from (30) that $\mathbb{E}(W_k^{(n)})^4 \leq C n^{2+2\delta}$ for all $k \geq 1$.

Also

$$\mathrm{Var}(W_1^{(n)}) = \mathbb{E}\left(\sum_{j=1}^{\lfloor n^\delta \rfloor} Y_j^{(n)}\right)^2$$

$$= \lfloor n^\delta \rfloor \mathbb{E}(Y_1^{(n)})^2 + 2 \sum_{j=1}^{\lfloor n^\delta \rfloor - 1} \operatorname{Cov}(Y_j^{(n)}, Y_{j+1}^{(n)}) \tag{32}$$

$$= \lfloor n^\delta \rfloor \mathbb{E}(Y_1^{(n)})^2 + 2(\lfloor n^\delta \rfloor - 1) \operatorname{Cov}(Y_1^{(n)}, Y_2^{(n)}).$$

In the above expression,

$$\begin{aligned} \mathbb{E}(Y_1^{(n)})^2 &= n \operatorname{Var}(Y_{1,1}) + 2 \sum_{s=1}^{n-1} \sum_{t=1}^{n-s} \operatorname{Cov}(Y_{s,1}, Y_{s+t,1}) \\ &= n \operatorname{Var}(Y_{1,1}) + 2 \sum_{s=1}^{n-1} (n-s) \operatorname{Cov}(Y_{1,1}, Y_{1+s,1}) \\ &= O(n) \qquad \text{as } n \to \infty, \end{aligned}$$

where the last equality follows because from the $\alpha$-mixing of the sequence $\{Y_{t,1}\}_{t\ge 1}$ we have $\sum_{t=2}^{\infty} \operatorname{Cov}(Y_{1,1}, Y_{t,1}) \le C \sum_{t=2}^{\infty} \alpha_t < \infty$ for some constant $C > 0$. Moreover, by the Cauchy–Schwarz inequality,

$$\operatorname{Cov}(Y_1^{(n)}, Y_2^{(n)}) \le \mathbb{E}(Y_1^{(n)})^2.$$

Thus, from (32), we have $\operatorname{Var}(W_1^{(n)}) = O(n^{1+\delta})$ as $n \to \infty$ and

$$\operatorname{Var}\left( \sum_{k=1}^{r_n} W_k^{(n)} \right) = O(n^{(1-\delta)+(1+\delta)}) = O(n^2) \qquad \text{as } n \to \infty. \tag{33}$$

Finally, for $0 < \delta < 1$,

$$\begin{aligned} &\lim_{n\to\infty} \sum_{k=1}^{r_n} \frac{1}{(\operatorname{Var} \sum_{k=1}^{r_n} W_k^{(n)})^2} \mathbb{E}(W_k^{(n)})^4 \\ &\quad \le \lim_{n\to\infty} \sum_{k=1}^{r_n} C \frac{n^{2+2\delta}}{n^4} = \lim_{n\to\infty} C n^{\delta-1} = 0. \end{aligned}$$

Thus by Lyapunov's central limit theorem [see Billingsley (1979), Theorem 27.3, page 312] we have that, for $0 < \delta < 1$, $1/(\sqrt{\sum_{k=1}^{r_n} \operatorname{Var}(W_k^{(n)})}) \sum_{k=1}^{r_n} W_k^{(n)}$ converges in probability to a standard normal random variable.

Now let $\eta_n := \sum_{k=1}^{r_n} \eta_k^{(n)}$. We will show that

$$\eta_n / n \to 0 \quad \text{in probability as} \quad n \to \infty. \tag{34}$$

Indeed,

$$\mathbb{E}(\eta_k^{(n)})^2 \le \sum_{i=1}^{n} \operatorname{Var}(Y_{i,k}) + 2n \sum_{i=2}^{n} \operatorname{Cov}(Y_{1,k}, Y_{i,k})$$

$$\leq n\,\mathrm{Var}(Y_{1,1}) + 2n\sum_{i=2}^{\infty} C_1 \exp(-C_2 i)$$

$$\leq Mn \qquad \text{for some constant } M > 0.$$

Thus, using the fact that $r_n = O(n^{1-\delta})$ as $n \to \infty$, we have, for $\varepsilon > 0$,

$$\mathbb{P}(|\eta_n| > n\varepsilon) \leq \frac{\mathbb{E}(\eta_n^2)}{n^2\varepsilon^2} = \frac{MnO(n^{1-\delta})}{n^2\varepsilon^2} \to 0 \qquad \text{as } n \to \infty.$$

This proves (34).

To complete the proof, we have to show that $\frac{E_n}{n} \to 0$ in probability as $n \to \infty$. First observe that number of terms in $E_n$ is at most $\lfloor n^{\delta} \rfloor$. Therefore taking $\delta = 1/2$, from (30) we have $\mathbb{E}(E_n^4) \leq Cn^3$. Hence, for $\varepsilon > 0$,

$$\mathbb{P}(|E_n| > n\varepsilon) \leq \frac{\mathbb{E}(E_n^4)}{n^4\varepsilon^4} \to 0 \qquad \text{as } n \to \infty. \tag{35}$$

Theorem 2.3(a) now follows by combining equations (34) and (35) and the fact that $\sum_{k=1}^{r_n} W_k{}^{(n)}/n$ has asymptotically a $N(0, s^2)$ distribution, where

$$s^2 = \mathrm{Var}(Y_{1,1}) + 2\sum_{i=2}^{\infty} \mathrm{Cov}(Y_{1,1}, Y_{i,1})$$

$$+ 2\sum_{i=1}^{\infty} \mathrm{Cov}(Y_{1,1}, Y_{i,2}) + 2\sum_{i=2}^{\infty} \mathrm{Cov}(Y_{1,2}, Y_{i,1}).$$

Note that to compute $s^2$ we use the fact that $\{(Y_{i,j}, Y_{i,j+1})\}_{i\geq 1}$ is an $\alpha$-mixing sequence.

## 6. Degree of a vertex.

To prove Proposition 2.1, observe that, given the vertex $(0,-1)$ is open, let

$$Y = \begin{cases} 1, & \text{if the vertex } (0,0) \text{ is open,} \\ 0, & \text{otherwise,} \end{cases}$$

$$X_1 = \#\{(i,0) : i \leq -1 : (i,0) \text{ is connected by an edge to } (0,-1)\},$$

$$X_2 = \#\{(i,0) : i \geq 1 : (i,0) \text{ is connected by an edge to } (0,-1)\}.$$

Clearly the degree of $(0,-1)$ equals $Y + X_1 + X_2$. Now given the vertex $(0,-1)$ is open, the probability that the vertex $(-l,0)$ is connected to $(0,-1)$ and that there are exactly $r-1$ vertices in $\{(i,0) : -l+1 \leq i \leq -1\}$ which are connected to $(0,-1)$ equals $\binom{l-1}{r-1} p^r (1-p)^{l-r} (1-p)^{2l-1}((1-p) + \frac{1}{2}p)$. Thus $\mathbb{P}(X_1 \geq r) = \sum_{l=r}^{\infty} \binom{l-1}{r-1} p^r (1-p)^{l-r} (1-p)^{2l-1}((1-p) + \frac{1}{2}p)$. An easy calculation now completes the proof of the proposition.

Similarly, in two dimensions, given that a vertex $v$ is open, the distribution of the number of edges of length $l$ "going up" from $v$ is binomial with parameters 2 and $(1-\frac{p}{2})(1-p)^{2l-1}$.

REMARK 6.1. From the above distributions we may calculate the quantities $\mathbb{E}(S_n)$, $\mathrm{Var}(S_n)$, $s^2$ and the related quantities involving $L_n$ required in Theorem 2.3 for two dimensions.

**Acknowledgments.** We are grateful to Professor S. Popov for his suggestions regarding the proof for $d = 3$ and to an anonymous referee for suggesting that the Burton–Keane argument would yield Theorem 2.2.

S. GANGOPADHYAY
INDIAN STATISTICAL INSTITUTE
STAT MATH DIVISION
203 B. T. ROAD
KOLKATA 700108
INDIA
E-MAIL: res9616@isical.ac.in

R. ROY
A. SARKAR
INDIAN STATISTICAL INSTITUTE
STAT MATH UNIT
7 S. J. S. SANSANWAL MARG
NEW DELHI 110016
INDIA
E-MAIL: rahul@isid.ac.in
E-MAIL: anish@isid.ac.in